# DEBUNKING CANTOR: NEW SET-THEORETICAL AND LOGICAL CONSIDERATIONS

JUAN A. PEREZ

ABSTRACT. For more than a century, Cantor's theory of transfinite numbers has played a pivotal role in set theory, with ramifications that extend to many areas of mathematics. This article extends earlier findings with a fresh look at the critical facts of Cantor's theory:

– Cantor's widely renowned Diagonalization Argument (CDA) is fully refuted by a set of counter-examples that expose the fallacy of this proof.

– The logical inconsistencies of CDA are revisited, exposing the short-comings of CDA's implementation of the *reductio* method of proof.

– The denumerability of the power set of the set of the natural numbers, $\mathcal{P}(\mathbb{N})$, is substantiated by a proof that takes full account of all the infinite subsets of $\mathbb{N}$. Such a result confirms the denumerability of the set of the real numbers, $\mathbb{R}$, and with it the countable nature of the continuum.

– Given that the denumerable character of (probably) all infinite sets makes their comparison in terms of one-to-one correspondences a rather pointless exercise, a new concept of relative cardinality is introduced which facilitates a quantitative evaluation of their different magnitudes.

## 1. INTRODUCTION

A previous report [15] presented a detailed and critical evaluation of the various proofs that underpin Cantor's theory of transfinite numbers [7,12,16]. Cantor's famous Diagonalization Argument (CDA) was particularly signalled for analysis, alongside other proofs supporting the uncountable nature of the set of real numbers, $\mathbb{R}$, and the power set of the set of natural numbers, $\mathcal{P}(\mathbb{N})$. Those proofs underpin much of modern set theory, with far reaching implications for most branches of mathematics. Consequently, their refutation (if correct) can be considered sufficiently important to merit further investigation. This article does precisely that, with a fresh look at the shortcomings of CDA, for which a number of counter-examples are described. The logical inadequacies of CDA are re-examined, reinforcing the previous analysis [15]. Furthermore, in order to confirm the denumerability of the power set of $\mathbb{N}$, $\mathcal{P}(\mathbb{N})$ (for which as many as three different proofs were already reported [15]), a new proof is described which takes into clear account all the infinite subsets of $\mathbb{N}$.

Since the filing of the original report, two other independent articles have reached the same conclusions, based on rationals that have much in common with our preceeding findings [4,10]. Our hope is that the new results presented here will further cement the inescapable conclusion that Cantorian mathematics needs to be expunged from the fabric of mathematical theory. The many implications for set theory and mathematical logic were extensively analysed before [15], so the interested reader is referred to the original material.

## 2. COUNTER-EXAMPLES OF CANTOR'S DIAGONALIZATION ARGUMENT

Cantor's Diagonalization Argument (CDA) [3,5,7,16] sits at the heart of his whole construction of transfinite number theory. Over the years, the simplicity of this argument has made it a favourite of set theorists and logicians alike [17], so it has been adapted to a great number of proofs. Hence, a refutation of CDA cannot be taken lightly. In order to analyse it in some detail in this and the following section, CDA will be reproduced here, adapted for the set of infinite binary strings [3,5,7]:

**Theorem.** The set of infinite binary strings is uncountable.

*Proof.* Suppose that the set $B$ of infinite binary strings is countable. Then we can list all the strings $s_n$ in $B$ as

$$s_1, s_2, s_3, \cdots, s_n, \cdots$$

with each string in $B$ appearing as $s_n$ for exactly one $n \in \mathbb{N}$, $n \geq 1$. We shall represent each string $s_n$ as

$$s_n = a_{n,1} a_{n,2} a_{n,3} \cdots a_{n,n} \cdots \quad n \in \mathbb{N}, n \geq 1$$

where each $a_{n,n}$ takes the value "0" or "1". We can then picture the set of strings $s_n$ written out in an array:

$$s_1 = a_{1,1} a_{1,2} a_{1,3} \cdots a_{1,n} \cdots$$
$$s_2 = a_{2,1} a_{2,2} a_{2,3} \cdots a_{2,n} \cdots$$
$$s_3 = a_{3,1} a_{3,2} a_{3,3} \cdots a_{3,n} \cdots$$
$$\vdots$$
$$s_n = a_{n,1} a_{n,2} a_{n,3} \cdots a_{n,n} \cdots$$
$$\vdots$$

Now define an "antidiagonal" string $s_{AD} = d_1 d_2 d_3 \cdots d_n \cdots$ by

$$d_n = \begin{cases} 1, \text{ if } a_{n,n} = 0, \\ 0, \text{ if } a_{n,n} = 1. \end{cases}$$

Then $s_D$ belongs to $B$. However, $s_{AD}$ has been constructed to disagree with each $s_n$ at the nth decimal place, so it cannot equal $s_n$ for any $n$. Thus $s_{AD}$ does not appear in the list, contradicting that the list contains all binary strings.

Therefore, we have that $B$ is uncountable. Q.E.D.

The main criticism originally raised against CDA was that the diagonal string $s_D$ can never "cover" the whole of the array [15]. This is best illustrated with a finite example: consider the set $B_4$ of all finite binary arrays of length 4, i.e. $s_n = a_{n,1}\, a_{n,2}\, a_{n,3}\, a_{n,4}$. It is simple to observe that the whole array consists of $2^4 = 16$ strings $s_n$, all of length 4:

| $s_1$ | $s_2$ | $s_3$ | $s_4$ | $s_5$ | $s_6$ | $s_7$ | $s_8$ | $s_9$ | $s_{10}$ | $s_{11}$ | $s_{12}$ | $s_{13}$ | $s_{14}$ | $s_{15}$ | $s_{16}$ |
|---|---|---|---|---|---|---|---|---|---|---|---|---|---|---|---|
| 0 | 1 | 0 | 0 | 0 | 1 | 1 | 1 | 0 | 0 | 0 | 1 | 1 | 1 | 0 | 1 |
| 0 | 0 | 1 | 0 | 0 | 1 | 0 | 0 | 1 | 1 | 0 | 1 | 1 | 0 | 1 | 1 |
| 0 | 0 | 0 | 1 | 0 | 0 | 1 | 0 | 1 | 0 | 1 | 1 | 0 | 1 | 1 | 1 |
| 0 | 0 | 0 | 0 | 1 | 0 | 0 | 1 | 0 | 1 | 1 | 0 | 1 | 1 | 1 | 1 |

It can be seen that the antidiagonal string will be $s_{AD} = 1\,1\,1\,1$, that is $s_{16}$, already a member of the array. Therefore, the constructed antidiagonal string only covers strings $s_1$, $s_2$, $s_3$ and $s_4$, unable to account for strings $s_5$ to $s_{16}$. In [15] a diagonal cover ($Dc$) was defined as the ratio between the sum of the members of the array covered by the antidiagonal string $s_{AD}$, divided by the total of members in the array. In the case of $B_4$, $Dc = 4/2^4 = 4/16$, significantly less than 1.

If the length of the binary strings were increased, the diagonal cover $Dc$ for a set $B_n$ of strings of length $n$ would be $Dc = n/2^n$ — therefore, the greater the value of $n$, the smaller the value of $Dc$, hence $Dc \ll 1$. For CDA to be correct when applied to the set of infinite binary strings $B$, it will be essential for $Dc = 1$. But this is evidently not the case, since

(2.1) $$Dc|_B = \lim_{n \to \infty} (Dc|_{B_n}) = \lim_{n \to \infty} (n/2^n) = 0$$

Perhaps the best way of refuting CDA would be to produce counter-examples that show conclusively that the antidiagonal string $s_D$ cannot cover the whole of the array of infinite binary strings. In view of the analysis presented here, finding counter-examples should be relatively simple. That this is the case can be shown by generating such counter-examples.

2.1. *Counter-examples using infinite binary strings.* If it is assumed that the set $B$ of infinite binary strings is denumerable, then it will be possible to construct an array following any ordering [7,13]:

(2.2)

|       | $s_1$ | $s_2$ | $s_3$ | $s_4$ | $s_5$ | $s_6$ | $s_7$ | $s_8$ | $s_9$ | $s_{10}$ | ... |
|-------|-------|-------|-------|-------|-------|-------|-------|-------|-------|----------|-----|
|       | **0** | 1 | 1 | 1 | 1 | 1 | 1 | 1 | 1 | 1 | ... |
|       | 0 | **0** | 1 | 1 | 1 | 1 | 1 | 1 | 1 | 1 | ... |
|       | 0 | 0 | **0** | 1 | 1 | 1 | 1 | 1 | 1 | 1 | ... |
|       | 0 | 0 | 0 | **0** | 1 | 1 | 1 | 1 | 1 | 1 | ... |
|       | 0 | 0 | 0 | 0 | **0** | 1 | 1 | 1 | 1 | 1 | ... |
|       | 0 | 0 | 0 | 0 | 0 | **0** | 1 | 1 | 1 | 1 | ... |
|       | 0 | 0 | 0 | 0 | 0 | 0 | **0** | 1 | 1 | 1 | ... |
|       | 0 | 0 | 0 | 0 | 0 | 0 | 0 | **0** | 1 | 1 | ... |
|       | 0 | 0 | 0 | 0 | 0 | 0 | 0 | 0 | **0** | 1 | ... |

It can be seen that, in this case, the antidiagonal string is $s_{AD} = 1111\cdots$, and this is the same infinite string that the strings in array $B$ tend towards, as the enumeration progresses undisturbed[1]. It it obvious that the antidiagonal string $s_{AD}$ will always be one step short of covering a string with the same number of 1s.

Another counter-example can be constructed, by changing in (2.2) all 1s for 0s, and all 0s for 1s:

(2.3)

|       | $s_1$ | $s_2$ | $s_3$ | $s_4$ | $s_5$ | $s_6$ | $s_7$ | $s_8$ | $s_9$ | $s_{10}$ | ... |
|-------|-------|-------|-------|-------|-------|-------|-------|-------|-------|----------|-----|
|       | **1** | 0 | 0 | 0 | 0 | 0 | 0 | 0 | 0 | 0 | ... |
|       | 1 | **1** | 0 | 0 | 0 | 0 | 0 | 0 | 0 | 0 | ... |
|       | 1 | 1 | **1** | 0 | 0 | 0 | 0 | 0 | 0 | 0 | ... |
|       | 1 | 1 | 1 | **1** | 0 | 0 | 0 | 0 | 0 | 0 | ... |
|       | 1 | 1 | 1 | 1 | **1** | 0 | 0 | 0 | 0 | 0 | ... |
|       | 1 | 1 | 1 | 1 | 1 | **1** | 0 | 0 | 0 | 0 | ... |
|       | 1 | 1 | 1 | 1 | 1 | 1 | **1** | 0 | 0 | 0 | ... |
|       | 1 | 1 | 1 | 1 | 1 | 1 | 1 | **1** | 0 | 0 | ... |
|       | 1 | 1 | 1 | 1 | 1 | 1 | 1 | 1 | **1** | 0 | ... |

On this occasion, the antidiagonal string is $s_{AD} = 0000\cdots$, and once again this is the same infinite string that the strings in array $B$ tend towards.

---

[1] This counter-example was independently recorded in an article [13] that only came to our attention when drafting this report. However, [13] concluded that no actual infinities can be considered in mathematical theory, and this is not a view we support.

The counter-examples constructed in (2.2) and (2.3) are not the only ones which can be conceived, an alternative construction will alternate 0s and 1s in $s_{AD}$:

(2.4)

| | $s_1$ | $s_2$ | $s_3$ | $s_4$ | $s_5$ | $s_6$ | $s_7$ | $s_8$ | $s_9$ | $s_{10}$ | $\cdots$ |
|---|---|---|---|---|---|---|---|---|---|---|---|
| | **0** | 1 | 1 | 1 | 1 | 1 | 1 | 1 | 1 | 1 | $\cdots$ |
| | 0 | **1** | 0 | 0 | 0 | 0 | 0 | 0 | 0 | 0 | $\cdots$ |
| | 0 | 0 | **0** | 1 | 1 | 1 | 1 | 1 | 1 | 1 | $\cdots$ |
| | 0 | 0 | 0 | **1** | 0 | 0 | 0 | 0 | 0 | 0 | $\cdots$ |
| | 0 | 0 | 0 | 0 | **0** | 1 | 1 | 1 | 1 | 1 | $\cdots$ |
| | 0 | 0 | 0 | 0 | 0 | **1** | 0 | 0 | 0 | 0 | $\cdots$ |
| | 0 | 0 | 0 | 0 | 0 | 0 | **0** | 1 | 1 | 1 | $\cdots$ |
| | 0 | 0 | 0 | 0 | 0 | 0 | 0 | **1** | 0 | 0 | $\cdots$ |
| | 0 | 0 | 0 | 0 | 0 | 0 | 0 | 0 | **0** | 1 | $\cdots$ |

In (2.4), the antidiagonal string is $s_{AD} = 1010101\cdots$, the same infinite string that the strings in array $B$ tend towards. And an alternative to (2.4) could be constructed by changing all the 1s for 0s and all the 0s for 1s:

(2.5)

| | $s_1$ | $s_2$ | $s_3$ | $s_4$ | $s_5$ | $s_6$ | $s_7$ | $s_8$ | $s_9$ | $s_{10}$ | $\cdots$ |
|---|---|---|---|---|---|---|---|---|---|---|---|
| | **1** | 0 | 0 | 0 | 0 | 0 | 0 | 0 | 0 | 0 | $\cdots$ |
| | 0 | **0** | 1 | 1 | 1 | 1 | 1 | 1 | 1 | 1 | $\cdots$ |
| | 0 | 0 | **1** | 0 | 0 | 0 | 0 | 0 | 0 | 0 | $\cdots$ |
| | 0 | 0 | 0 | **0** | 1 | 1 | 1 | 1 | 1 | 1 | $\cdots$ |
| | 0 | 0 | 0 | 0 | **1** | 0 | 0 | 0 | 0 | 0 | $\cdots$ |
| | 0 | 0 | 0 | 0 | 0 | **0** | 1 | 1 | 1 | 1 | $\cdots$ |
| | 0 | 0 | 0 | 0 | 0 | 0 | **1** | 0 | 0 | 0 | $\cdots$ |
| | 0 | 0 | 0 | 0 | 0 | 0 | 0 | **0** | 1 | 1 | $\cdots$ |
| | 0 | 0 | 0 | 0 | 0 | 0 | 0 | 0 | **1** | 0 | $\cdots$ |

In (2.5), the antidiagonal string will be $s_{AD} = 0101010\cdots$, once more the same infinite string that the strings in array $B$ tend towards.

It should be equally obvious that, in all constructions (2.2) to (2.5), it will be possible to place a completely randomised set of 1s and 0s below the diagonal line, leaving the antidiagonal string $s_{AD}$ unaffected, thus highlighting that the number of possible counter-examples is, in fact, infinite. For example, illustrating this point, an alternative to (2.2) might be:

(2.6)

|     | $s_1$ | $s_2$ | $s_3$ | $s_4$ | $s_5$ | $s_6$ | $s_7$ | $s_8$ | $s_9$ | $s_{10}$ | $\cdots$ |
|---|---|---|---|---|---|---|---|---|---|---|---|
| | 0. | 1 | 1 | 1 | 1 | 1 | 1 | 1 | 1 | 1 | $\cdots$ |
| | 0 | 0. | 1 | 1 | 1 | 1 | 1 | 1 | 1 | 1 | $\cdots$ |
| | 1 | 0 | 0. | 1 | 1 | 1 | 1 | 1 | 1 | 1 | $\cdots$ |
| | 1 | 1 | 0 | 0. | 1 | 1 | 1 | 1 | 1 | 1 | $\cdots$ |
| | 1 | 0 | 1 | 1 | 0. | 1 | 1 | 1 | 1 | 1 | $\cdots$ |
| | 0 | 1 | 1 | 0 | 1 | 0. | 1 | 1 | 1 | 1 | $\cdots$ |
| | 0 | 1 | 0 | 1 | 0 | 0 | 0. | 1 | 1 | 1 | $\cdots$ |
| | 1 | 0 | 0 | 1 | 1 | 1 | 0 | 0. | 1 | 1 | $\cdots$ |
| | 0 | 1 | 1 | 0 | 1 | 0 | 1 | 0 | 0. | 1 | $\cdots$ |

originating the same antidiagonal string as in (2.2), i.e. $s_{AD} = 1\,1\,1\,1\cdots$. And the same can be said of arrays (2.3) to (2.5).

2.2. *Counter-examples using decimal representations of the reals.* It is usual to apply CDA to sets of decimal representations of real numbers [7,12], so it makes sense to construct counter-examples for this scenario. Consider $R$, the set of real numbers in the interval [0,1), and assume it is denumerable. Accordingly, it will be possible to write the array:

(2.7) $$R = \{r_1, r_2, r_3, r_4, ..., r_n, ...\}$$

where each real in [0,1) appears as $r_n$ for exactly one $n \in \mathbb{N}$, $n \geq 1$. Each real $r_n$ can be represented by its decimal expansion, such that

(2.8) 
$$r_1 = 0.a_{1,1}\,a_{1,2}\,a_{1,3}\,a_{1,4}\cdots a_{1,n}\cdots$$
$$r_2 = 0.a_{2,1}\,a_{2,2}\,a_{2,3}\,a_{2,4}\cdots a_{2,n}\cdots$$
$$r_3 = 0.a_{3,1}\,a_{2,2}\,a_{3,3}\,a_{3,4}\cdots a_{3,n}\cdots$$
$$r_2 = 0.a_{4,1}\,a_{4,2}\,a_{4,3}\,a_{4,4}\cdots a_{4,n}\cdots$$
$$\vdots$$
$$r_n = 0.a_{n,1}\,a_{n,2}\,a_{n,3}\,a_{n,4}\cdots a_{n,n}\cdots$$
$$\vdots$$

where the digits $a_{n,n}$ take any of the values from 0 to 9, but avoiding the use of recurring 9s [6] (so that e.g. 0.2 is represented by 0.2000..., rather than 0.1999...). Since $R$ is denumerable, the array (2.8) can be listed in any order of our choice [7,13]; therefore, it will be acceptable to initiate the ordering as follows:

(2.9)
$$\begin{aligned}
r_0 &= 0.\,0\ 0\ 0\ 0\ 0\ 0\ 0\ 0\ 0\ 0\ \cdots \\
r_1 &= 0.\,a_1\ 0\ 0\ 0\ 0\ 0\ 0\ 0\ 0\ 0\ \cdots \\
r_2 &= 0.\,a_1\ a_2\ 0\ 0\ 0\ 0\ 0\ 0\ 0\ 0\ \cdots \\
r_3 &= 0.\,a_1\ a_2\ a_3\ 0\ 0\ 0\ 0\ 0\ 0\ 0\ \cdots \\
r_4 &= 0.\,a_1\ a_2\ a_3\ a_4\ 0\ 0\ 0\ 0\ 0\ 0\ \cdots \\
r_5 &= 0.\,a_1\ a_2\ a_3\ a_4\ a_5\ 0\ 0\ 0\ 0\ 0\ \cdots \\
r_6 &= 0.\,a_1\ a_2\ a_3\ a_4\ a_5\ a_6\ 0\ 0\ 0\ 0\ \cdots \\
r_7 &= 0.\,a_1\ a_2\ a_3\ a_4\ a_5\ a_6\ a_7\ 0\ 0\ 0\ \cdots \\
r_8 &= 0.\,a_1\ a_2\ a_3\ a_4\ a_5\ a_6\ a_7\ a_8\ 0\ 0\ \cdots \\
r_9 &= 0.\,a_1\ a_2\ a_3\ a_4\ a_5\ a_6\ a_7\ a_8\ a_9\ 0\ \cdots \\
r_{10} &= 0.\,a_1\ a_2\ a_3\ a_4\ a_5\ a_6\ a_7\ a_8\ a_9\ a_{10}\ \cdots \\
&\ \ \vdots
\end{aligned}$$

where the digits $a_n$ take values from 1 to 9, to be specified by the choices made for the digits of the antidiagonal number $r_{AD} = 0.\,a_1 a_2 a_3 a_4 a_5 a_6 a_7 a_8 a_9 a_{10}\ldots$ In other words, the enumeration of the array $R$ is written based on the selections of digits $a_n$ for $r_{AD}$. It is easy to see that, as the construction of both the array $R$ and the antidiagonal number $r_{AD}$ progresses unimpeded, $r_{AD}$ will be incorporated in $R$, contradicting CDA. It can also be observed that the construction (2.9) originates an infinite number of counter-examples, since the digits $a_n$ can take any combination of values 1 to 9, provided that $a_{n,n} \neq 0$. Furthermore, the 0s below the diagonal could be replaced by any random combination of values 0 to 9 (while leaving $r_{AD}$ unaltered), hence increasing even more the myriad of counter-examples that (2.9) provides.

The counter-examples (2.2) to (2.6), and (2.9), do more than enough to fully refute CDA.

## 3. LOGICAL SHORT-COMINGS OF CANTOR'S DIAGONALIZATION ARGUMENT

The refutation of CDA suggests that something must have been adrift with its implementation as a proof (by contradiction). This much was evaluated in [15]. However, an inaccuracy was made in the original analysis that left matters in an unsatisfactory state, hence a further evaluation is warranted. Although many of the original conclusions still remained valid, a more nuanced approach was required, which is undertaken here.

### 3.1. *Proofs by contradiction (reductio ad absurdum).*
In an attempt to prove a given statement $P$, a proof by contradiction essentially works by first assuming the truth of the opposite statement, this is the negation of $P$ ($\neg P$), and then allowing the implementation of standard rules of inference to proceed through a string of interconnecting statements $Q_1$, $Q_2$, $\cdots$, $Q_n$ until a final contradiction is reached (i.e. a statement which is always false) [8,14]. Since a contradiction is a statement than can never be true (it is commonplace to describe it as a composite statement, $R \wedge \neg R$, hence reinforcing its falsehood), its negation leads to the negation of $\neg P$, and this, in turn, to the truth of $P$ [8]. The associated chain of inference can be written as:

$$\neg P \Rightarrow Q_1 \Rightarrow Q_2 \Rightarrow \cdots \Rightarrow Q_n \Rightarrow (R \wedge \neg R) \tag{3.1}$$

so the rule of hypothegical syllogism [14] implies

$$\neg P \Rightarrow (R \wedge \neg R) \tag{3.2}$$

and, by *modus tollens* and double negation [14],

$$\neg (R \wedge \neg R) \Rightarrow \neg(\neg P) \Rightarrow P \tag{3.3}$$

completing the proof. A variation on this theme reported in [15] has, as the final statement in the chain of inference, the initial statement $P$, this is

$$\neg P \Rightarrow Q_1 \Rightarrow Q_2 \Rightarrow \cdots \Rightarrow Q_n \Rightarrow P \tag{3.4}$$

so the rule of hypothegical syllogism combined with conjunction introduction [14] now implies that

$$\neg P \Rightarrow (P \wedge \neg P) \tag{3.5}$$

and, once more by *modus tollens* and double negation [14],

$$\neg (P \wedge \neg P) \Rightarrow \neg(\neg P) \Rightarrow P \tag{3.6}$$

and the proof is again complete. The chain of inference (3.4) is relevant to our analysis, given that this is the form of proof associated with CDA [15].

One fundamental aspect of proofs by contradiction is the fact that, in order to derive the truth of $P$, the truth of all the intermediate statements $Q_n$ in (3.1), or (3.4), has to be independently asserted. Quoting from [8]:

> "Such a proof (*reductio ad absurdum*) consists of a deduction of a contradiction from the negation of the statement whose proof is required. That this is a legitimate procedure (..) can be seen as follows. If we have an argument which is known to be an instance of a valid argument form, and its conclusion is known to be false, then at least one of the premises must be false, If all the premises are known to be true except one (the assumed one), then the legitimate deduction is that this assumed one is the one which is false."

Such a prerequisite is fundamental to the success of these proofs. The chains of inference (3.1) or (3.4) do not offer any additional complication, but the same cannot be said of chains of inference where the connectors are biconditional instead of single conditional:

(3.7) $\quad\neg P \Leftrightarrow Q_1 \Leftrightarrow Q_2 \Leftrightarrow \cdots \Leftrightarrow Q_n \Rightarrow (R \wedge \neg R)$

(3.8) $\quad\neg P \Leftrightarrow Q_1 \Leftrightarrow Q_2 \Leftrightarrow \cdots \Leftrightarrow Q_n \Rightarrow P$

In (3.7) and (3.8), the truth of the intermediate statements $Q_1,\ldots,Q_n$ is directly associated to the truth of $\neg P$, since they all are equivalent statements [8]. Therefore, if $\neg P$ is false, so are $Q_1,\ldots,Q_n$. In other words, the falsehood of $(R \wedge \neg R)$, or $(P \wedge \neg P)$, implies the falsehood of them all, $Q_1,\ldots,Q_n$ as well as $\neg P$. Consequently, the proofs fail to have a single true statement underpinning the sought conclusion, i.e. the truth of $P$. It is hard to see this scenario being nothing but a corruption of the method of proof by contradiction. This much was concluded in [15].

However, there is a "half-way house" situation where having biconditional statements connecting $\neg P$ to some, but not all, of the statements $Q_1,\ldots,Q_n$ does not compromise the validity of the proof:

(3.9) $\quad\neg P \Leftrightarrow Q_1 \Leftrightarrow Q_2 \Leftrightarrow \cdots \Leftrightarrow Q_{i-1} \Leftrightarrow Q_i \Rightarrow Q_{i+1} \Rightarrow \cdots \Rightarrow Q_n \Rightarrow (R \wedge \neg R)$

(3.10) $\quad\neg P \Leftrightarrow Q_1 \Leftrightarrow Q_2 \Leftrightarrow \cdots \Leftrightarrow Q_{i-1} \Leftrightarrow Q_i \Rightarrow Q_{i+1} \Rightarrow \cdots \Rightarrow Q_n \Rightarrow P$

In (3.9) and (3.10), the truth of the statements $Q_{i+1},\ldots,Q_n$ is not associated to the truth of $\neg P$ (unlike $Q_1,\ldots,Q_i$), and that will be sufficient to validate the proof, provided the statements $Q_{i+1},\ldots,Q_n$ were shown to be true. It is this observation what we failed to notice in our original report [15][2].

Knowing already that CDA is a flawed proof, we are now in a position to evaluate the logical structure of CDA. If we take the presentation of CDA already described in Section 2, we could dissect the chain of inference as follows:

- $P$ = 'The set $B$ of infinite binary strings is uncountable'
- $\neg P$ = 'The set $B$ of infinite binary strings is countable'
- $Q_1$ = 'The strings $s_n$ in $B$ can be listed as
  $$s_1, s_2, s_3, \ldots, s_n, \ldots$$
  where $n \in \mathbb{N}, n \geq 1$'
- $Q_2$ = 'We can picture the set of strings $s_n$ written out in an array:
  $$s_n = a_{n,1}\, a_{n,2}\, a_{n,3} \ldots a_{n,n} \ldots$$
  where $n \in \mathbb{N}, n \geq 1$'
- $Q_3$ = 'We define an "antidiagonal" string $s_{AD} = d_1\, d_2\, d_3 \ldots d_n \ldots$ by
  $$d_n = \begin{cases} 1, & \text{if } a_{n,n} = 0 \\ 0, & \text{if } a_{n,n} = 1 \end{cases}$$

---

[2] Fortunately, all the proofs of nondenumerability that were analysed in [15] fall into the category of (3.7) or (3.8), so the conclusions reported there remain sound.

$s_{AD}$ belongs to $B$ but is $s_{AD} \neq a_{n,n}$ for all $n \in \mathbb{N}$, so it cannot be part of the array'

- $Q_4$ = 'The array is not a complete listing of the elements of $B$'

The list of statements forms the logical sequence

(3.11) $$\neg P \Leftrightarrow Q_1 \Leftrightarrow Q_2 \Leftrightarrow Q_3 \Rightarrow Q_4 \Leftrightarrow P$$

where the connectives linking $\neg P$ with $Q_1$, $Q_2$ and $Q_3$ are all biconditional, leaving just a single conditional connective between $Q_3$ and $Q_4$ since, in principle, there could be other reasons (not addressed by the proof) why the array is not a complete listing of $B$. The final connective between $Q_4$ and $P$ is also biconditional. It is important to understand that the connective between statements $Q_2$ and $Q_3$ is biconditional: the antidiagonal string $s_{AD}$ can only be defined based on the construction of the array and, in reverse, the definition of $s_{AD}$ implies the existence of the countable array.

The chain of inference (3.11) is an example of (3.8), without a single true intermediate statement underpinning the validity of the proof. Since we already know that CDA is flawed, it should come as no suprise that its logical structure fails to meet the requirements of a correct proof by contradiction. In fact, this failure could have been used to point to the short-comings of CDA. Since it is the case that $Q_3 \Rightarrow (P \wedge \neg P)$, there are no circumstances under which $Q_3$ can be a true statement.

The implications of the flawed nature of CDA as a method of proof are considerable. Diagonalisation arguments have been used extensively by set-theorists and logicians over the years, and quite a number of important results (including Gödel's famous theorems of incompleteness [17]) are underpinned by such arguments. This issue was comprehensively analysed in our previous report [15].

With regard to the nondenumerability of the set of real numbers, $\mathbb{R}$, non-set-theoretical proofs can be found in the literature that come from other branches of mathematics [11]. It will be of interest to verify whether such proofs also lack a reliable logical structure.

In our previous report [15], we introduced the definition of *inconceivable* statements, to be used as a preventative measure against the construction of incorrect proofs:

**Definition 3.1.** *A mathematical statement $Q$ is said to be inconceivable when there is another statement $P$ such that*

***i)*** *$(Q \Rightarrow P) \wedge (Q \Rightarrow \neg P)$, or  **ii)** $Q \Rightarrow ((P \Rightarrow \neg P) \vee (\neg P \Rightarrow P))$*

*Otherwise, the statement $Q$ is considered conceivable.*

This definition lead to the formulation of a Principle (of *Conceivable Proof*) that needs a slight alteration, in order to account for proof constructions such as (3.9) and (3.10).

**Principle 5.2 (of Conceivable Proof).** *No mathematial proof by contradiction can be judged valid if (in the absence of any true statement or statements underpinning the proof) its construction includes one or more inconceivable statements; an exception will be when the purpose of the proof is to demonstrate the falsehood of such an inconceivable statement, provided that the resulting contradiction is not conceptually linked to the initial assumption of the proof.*

This principle was initially put forward to prevent the construction of erroneous proofs like CDA [15]. Proofs by contradiction are everywhere in the mathematical literature, thus the prevention of mistakes in their formulation seems warranted.

## 4. Denumerability of the Power Set of $\mathbb{N}$ ($\mathcal{P}(\mathbb{N})$)

The refutation of CDA reported here, as well as the preceding critical evaluations of the remaining set-theoretical proofs on the uncountability of $\mathbb{R}$ and $\mathcal{P}(\mathbb{N})$ [15], do not, just by themselves, prove that these sets are denumerable. Such a conclusion can only be reached with the construction of the relevant proof/s. To this effect, in our previous report we described three independent proofs of the denumerability of $\mathcal{P}(\mathbb{N})$ [15], whose conpletion required the formulation of a new theorem (of actual countable infinity) [15], a natural extension of the axiom of infinity [7,12,16]. A new proof is presented here that does not make use of such a theorem, and takes full account of all the infinite subsets of $\mathbb{N}$.

Two preliminary facts need to be taken into account before dealing with the proof, which will be used in its construction:

*i)* Firstly, consider a well-known theorem for the union of countable sets [2,16]:

> **Theorem 4.1.** *If $A_n$ is a countable set for each $n \in \mathbb{N}$, then the union $A := \bigcup_{n=1}^{\infty} A_n$ is countable.*

In Theorem 4.1, for countable it will also be understood infinitely countable, i.e. denumerable.

*ii)* Secondly, consider a finite set $A_n$ of $n$ members, with $n \in \mathbb{N}$, and also consider its power set, $\mathcal{P}(A_n)$, with cardinality given by [1,15]:

$$(4.1) \quad |\mathcal{P}(A_n)| = 2^n = \binom{n}{0} + \binom{n}{1} + \cdots + \binom{n}{p} + \cdots + \binom{n}{n/2} + \cdots + \binom{n}{n-p} + \cdots + \binom{n}{n-1} + \binom{n}{n}$$

where each binomial coefficient equals the cardinality of a given subset of $\mathcal{P}(A_n)$.[3] If we name as $\mathcal{N}(A_p)$ the subset of $\mathcal{P}(A_n)$ formed by all the subsets of $A_n$ with cardinality $p$, it will be the case that

$$(4.2) \qquad |\mathcal{N}(A_p)| = \binom{n}{p}$$

and we will be able to write that

$$(4.3) \qquad \mathcal{P}(A_n) = \bigcup_{p=0}^{n} \mathcal{N}(A_p)$$

and

$$(4.4) \qquad |\mathcal{P}(A_n)| = \sum_{p=0}^{n} |\mathcal{N}(A_p)| = \sum_{p=0}^{n} \binom{n}{p}$$

since the subsets $\mathcal{N}(A_p)$ are pairwise disjoint.

A property of the binomial coefficients, to be used in what follows, is the relationship between consecutive coefficients [1], that is

$$(4.5) \qquad \binom{n}{p+1} = \binom{n}{p} \cdot \frac{(n-p)}{(p+1)}$$

Applying (4.5) to consecutive binomial coefficients beyond $n/2$, it can be easily deduced that

$$(4.6) \qquad \binom{n}{n/2+d+1} = \binom{n}{n/2+d} \cdot q$$

where the ratio $q$ ($q \in \mathbb{Q}$) is given by

$$(4.7) \qquad q = \frac{(n-2d)}{[n+2(d+1)]} = \frac{1-(2/n) \cdot d}{[1+(2/n) \cdot (d+1)]}$$

with $0 < q < 1$, and $0 \leq d \leq n/2 - 1$ ($d \in \mathbb{N}$). The ratio $q$ will take the limit values:

$$(4.8) \quad d = 0 \;\rightarrow\; \binom{n}{n/2+1} = \binom{n}{n/2} \cdot q \;\text{ with }\; q = \frac{1}{1+(2/n)}$$

$$(4.9) \quad d = n/2-1 \;\rightarrow\; \binom{n}{n} = \binom{n}{n-1} \cdot q \;\text{ with }\; q = 1/n$$

The above results can be illustrated graphically with an example, e.g. $n = 40$ (Figure 1). It is also relevant to examine the values that $q$ takes for a range of values of $d$ as a function of $n$ (Table 1).

---

[3] In (4.1) it has been assumed, for simplicity, that $n$ is an even natural, i.e. $n \equiv 0 \pmod{2}$, so that the binomial expansion of $2^n$ has just one single central term $\binom{n}{n/2}$.

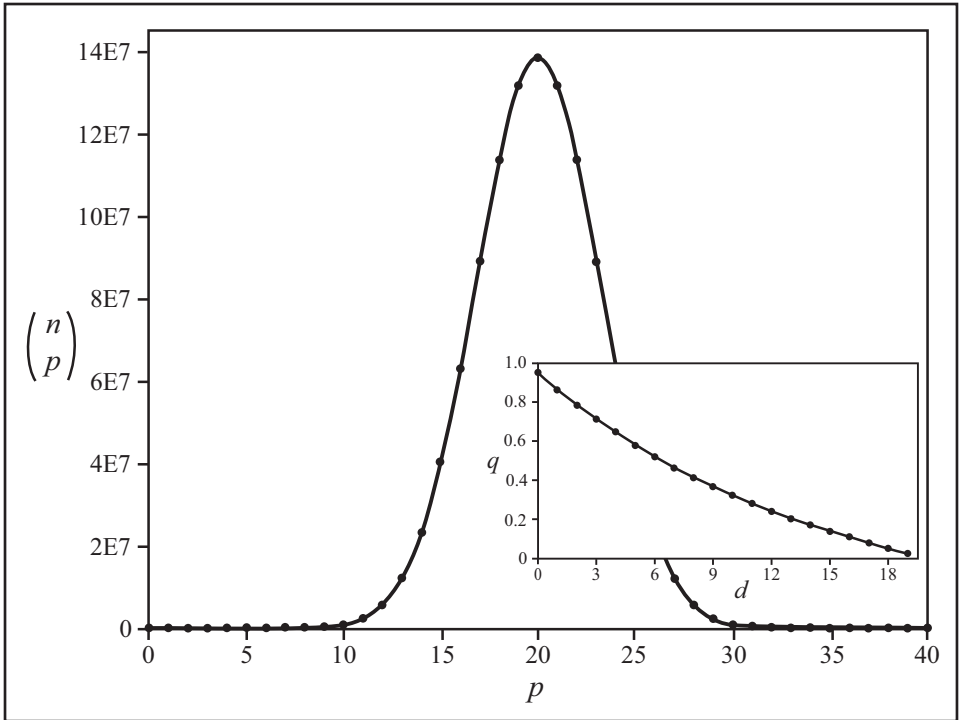

FIGURE 1. $\binom{n}{p}$ vs $p$ and $q$ vs $d$ for $n=40$.

| $d$ | $q$ | $d$ | $q$ | $d$ | $q$ |
|---|---|---|---|---|---|
| 0 | $1/(1+2/n)$ | $n/10$ | $2/(3+5/n)$ | $n/2-8$ | $8/(n-7)$ |
| 1 | $(1-2/n)/(1+4/n)$ | $n/9$ | $7/(11+18/n)$ | $n/2-7$ | $7/(n-6)$ |
| 2 | $(1-4/n)/(1+6/n)$ | $n/8$ | $3/(5+8/n)$ | $n/2-6$ | $6/(n-5)$ |
| 3 | $(1-6/n)/(1+8/n)$ | $n/7$ | $5/(9+14/n)$ | $n/2-5$ | $5/(n-4)$ |
| 4 | $(1-8/n)/(1+10/n)$ | $n/6$ | $2/(4+6/n)$ | $n/2-4$ | $4/(n-3)$ |
| 5 | $(1-10/n)/(1+12/n)$ | $n/5$ | $3/(7+10/n)$ | $n/2-3$ | $3/(n-2)$ |
| 6 | $(1-12/n)/(1+14/n)$ | $n/4$ | $1/(3+4/n)$ | $n/2-2$ | $2/(n-1)$ |
| 7 | $(1-14/n)/(1+16/n)$ | $n/3$ | $1/(5+6/n)$ | $n/2-1$ | $1/n$ |

TABLE 1. Values of $q$ vs $d$ as a function of $n$.

With these two facts now established, we are in a position to describe our new proof of the denumerability of the power set $\mathcal{P}(\mathbb{N})$.

**Theorem 4.2 (Denumerability of the Power Set)**. *Let $\mathbb{N}$ be the set of all natural numbers, $\mathbb{N} = \{0, 1, 2, 3, \cdots, p, \cdots\}$. Its power set $\mathcal{P}(\mathbb{N})$, the set of all subsets of $\mathbb{N}$, is denumerable*.

*Proof.* The power set $\mathcal{P}(\mathbb{N})$ will be the union of all sets $\mathcal{N}(N_p)$, i.e. the sets of all subsets of $\mathbb{N}$ with cardinality $p$:

$$(4.10) \qquad \mathcal{P}(\mathbb{N}) = \bigcup_{p=0}^{\aleph_0} \mathcal{N}(N_p)$$

where $\aleph_0$ denotes the cardinality of $\mathbb{N}$, i.e. $|\mathbb{N}| = \aleph_0$ [7,12,16]. Therefore, the sets $\mathcal{N}(N_n)$ will comprise both finite and infinite sets. Since all sets $\mathcal{N}(N_p)$ are pairwise disjoint, the cardinality of $\mathcal{P}(\mathbb{N})$ can be expressed as the summation of the cardinalities of all sets $\mathcal{N}(N_p)$:

$$(4.11) \qquad |\mathcal{P}(\mathbb{N})| = \sum_{p=0}^{\aleph_0} |\mathcal{N}(N_p)|$$

with the total count of sets $\mathcal{N}(N_p)$ being denumerable, i.e. if we denote as $\mathcal{T}_{\aleph_0}$ the set whose members are all the sets $\mathcal{N}(N_p)$:

$$(4.12) \qquad \mathcal{T}_{\aleph_0} = \{\mathcal{N}(N_0), \mathcal{N}(N_1), \mathcal{N}(N_2), \ldots, \mathcal{N}(N_p), \ldots, \mathcal{N}(N_{\aleph_0})\}$$

it will be the case that

$$(4.13) \qquad |\mathcal{T}_{\aleph_0}| = \aleph_0$$

In order to prove that the power set $\mathcal{P}(\mathbb{N})$ is denumerable, it will be sufficient to prove that each set $\mathcal{N}(N_p)$ is denumerable, that is

$$(4.14) \qquad |\mathcal{N}(N_p)| = \aleph_0 \; \forall p \in \mathbb{N} \; \wedge \; \forall \mathcal{N}(N_p) \colon |N_p| = \aleph_0$$

It is already well-documented that the set whose members are all the finite subsets of $\mathbb{N}$, denoted $\mathcal{F}(\mathbb{N})$, is denumerable [16]: A proof of this statement can easily be constructed using Theorem 4.1 and mathematical induction:

– $|\mathcal{N}(N_0)| = 1$, since the only member of $\mathcal{N}(N_0)$ is the empty set, $\emptyset$.
– $|\mathcal{N}(N_1)| = \aleph_0$, as the members of $\mathcal{N}(N_1)$ are all the singletons $\{i\}$, $\forall i \in \mathbb{N}$.
– $|\mathcal{N}(N_2)| = \aleph_0$. The members of $\mathcal{N}(N_2)$, i.e. the pairs $\{i, j\}$, $\forall i, j \in \mathbb{N}$, can be constructed by the union of all the singletons $\{i\}$ with the singletons $\{j\}$, provided that, to avoid repetition, $j > i$. For each singleton $\{i\}$, a denumerable set of pairs $\{i, j\}$ will be generated. Since the set of singletons $\{i\}$ is also denumerable, application of Theorem 4.1 will imply that the total of pairs $\{i,j\}$ is indeed denumerable.
– The inductive step: If it is assumed that the set $\mathcal{N}(N_p)$, $p \in \mathbb{N}$, is denumerable, this will imply that the set $\mathcal{N}(N_{p+1})$ is also denumerable. The members of $\mathcal{N}(N_p)$ will be all the subsets of $\mathbb{N}$ with cardinality $p$, e.g. $\{0,1,2,3,\cdots,p-1\}$.

To construct the members of $\mathcal{N}(N_{p+1})$, we will need to undertake the union of every member of $\mathcal{N}(N_p)$ (this is, a subset of $\mathbb{N}$ with cardinality $p$) with every singleton $\{i\}$, such that $i$ is greater than any of the elements of the given subset, in order to avoid repetitions. Therefore, for every member of $\mathcal{N}(N_p)$, a denumerable list of subsets of $\mathbb{N}$ with cardinality $p+1$ will be generated. Since it has been assumed that $\mathcal{N}(N_p)$ is a denumerable set, the application of Theorem 4.1 will conclude that $\mathcal{N}(N_{p+1})$ is also denumerable.

– A final application of Theorem 4.1 will help us to reach the conclusion that $\mathcal{F}(\mathbb{N})$ is denumerable.

All members of $\mathcal{F}(\mathbb{N})$ are finite subsets of $\mathbb{N}$. To construct the infinite subsets of $\mathbb{N}$, a starting point can be those subsets that are obtained by extracting the differences $\mathbb{N} \setminus N_p$, this is, extracting from $\mathbb{N}$ the finite subsets $N_p$. To illustrate this point, two examples of $\mathbb{N} \setminus N_1$ will be the subsets of $\mathbb{N}$ $\{1,2,3,4,\ldots,p,\ldots\}$ and $\{0,2,3,4,\ldots,p,\ldots\}$, while two examples of $\mathbb{N} \setminus N_2$ will be the subsets of $\mathbb{N}$ $\{2,3,4,\ldots,p,\ldots\}$ and $\{0,3,4,\ldots,p,\ldots\}$. And so on. Since a one-to-one correspondence can be established between every $\mathcal{N}(N_p)$ and $\mathcal{N}(N_{\mathbb{N}\setminus p})$ (where $N_{\mathbb{N}\setminus p}$ denotes the set of subsets of $\mathbb{N}$ obtained by extracting the differences $\mathbb{N}\setminus N_p$, for a given cardinality $p$), i.e. $\mathcal{N}(N_p) \leftrightarrow \mathcal{N}(N_{\mathbb{N}\setminus p})$, it will always be the case that

$$(4.15) \qquad \left|\mathcal{N}(N_p)\right| = \left|\mathcal{N}(N_{\mathbb{N}\setminus p})\right| = \aleph_0 \quad \forall p \in \mathbb{N} \land p \neq 0$$

Returning to (4.11), expand this statement as follows:

$$(4.16) \quad \left|\mathcal{P}(\mathbb{N})\right| = \left|\mathcal{N}(N_0)\right| + \left|\mathcal{N}(N_1)\right| + \cdots + \left|\mathcal{N}(N_p)\right| + \cdots + \left|\mathcal{N}(N_{\mathbb{E},\mathbb{O}})\right| + \\ + \cdots + \left|\mathcal{N}(N_{\mathbb{N}\setminus p})\right| + \cdots + \left|\mathcal{N}(N_{\mathbb{N}\setminus 1})\right| + \left|\mathcal{N}(N_{\mathbb{N}\setminus 0})\right|$$

where $N_{\mathbb{E},\mathbb{O}}$ denotes the set of infinite subsets of $\mathbb{N}$ with equal numbers of members of $\mathbb{N}$ missing as showing. Two examples of these subsets are the set of all even numbers, $\mathbb{E} = \{0,2,4,6,8,\ldots\}$, and the set of all odd numbers, $\mathbb{O} = \{1,3,5,7,9,\ldots\}$ (notice that $\mathbb{E} \cup \mathbb{O} = \mathbb{N}$, and $\mathbb{E} \cap \mathbb{O} = \emptyset$). We know already that some of the sets (of subsets of $\mathbb{N}$) are denumerable:

$$(4.17) \quad \left|\mathcal{P}(\mathbb{N})\right| = 1 + \aleph_0 + \cdots + \aleph_0 + \cdots + \left|\mathcal{N}(N_{\mathbb{E},\mathbb{O}})\right| + \cdots + \aleph_0 + \cdots + \aleph_0 + 1$$

Adapting (4.4) from the power set of a finite set, $A_n$, to the power set of $\mathbb{N}$, the cardinality of $\mathcal{P}(\mathbb{N})$ will be given by

$$(4.18) \qquad \left|\mathcal{P}(\mathbb{N})\right| = \lim_{n \to \aleph_0} \left[ \sum_{p=0}^{\aleph_0} \binom{n}{p} \right] = \sum_{p=0}^{\aleph_0} \left[ \lim_{n \to \aleph_0} \binom{n}{p} \right]$$

If (4.16) and (4.18) are compared, it becomes possible to derive the cardinality of $\mathcal{N}(N_{\mathbb{E},\mathbb{O}})$ as given by

$$\left|\mathcal{N}(N_{\mathbb{E},\mathbb{O}})\right| = \lim_{n \to \aleph_0} \binom{n}{n/2} \tag{4.19}$$

since $\left|\mathcal{N}(N_{\mathbb{E},\mathbb{O}})\right|$ is the central term of the expansion (4.16). For $\mathcal{P}(\mathbb{N})$ to be uncountable, it will be required that at least one of the members of the expansion (4.16) is uncountable. Accordingly, $\mathcal{N}(N_{\mathbb{E},\mathbb{O}})$ will have to be uncountable, as it is the central term of (4.16). If the arithmetics of transfinite cardinals are taken into consideration [7,12,16], it will be realised that (4.19) does not generate a transfinite cardinal larger than $\aleph_0$:

$$\left|\mathcal{N}(N_{\mathbb{E},\mathbb{O}})\right| = \lim_{n \to \aleph_0} \binom{n}{n/2} = \lim_{n \to \aleph_0} \left[\frac{n!}{(n/2)! \cdot (n/2)!}\right] = \aleph_0 \tag{4.20}$$

since $\aleph_0 \cdot \aleph_0 = \aleph_0$ [7,12,16]; as the limit is taken, the factorial $n!$ will never be able to grow in value beyond $\aleph_0$.

What is true of $\mathcal{N}(N_{\mathbb{E},\mathbb{O}})$, will also be true of all the subsets of $\mathcal{P}(\mathbb{N})$ whose members are infinite subsets of $\mathbb{N}$. Therefore, all of them will be denumerable. To corroborate this conclusion, we can proceed with the following analysis: Assume that $\mathcal{N}(N_{\mathbb{E},\mathbb{O}})$ is an uncountable set, so that its cardinality is $\aleph_1$, i.e. the least uncountable cardinal [7,12,16]. According to the arithmetics of transfinite cardinals, it is the case that $\aleph_1 \cdot k = \aleph_1$, $\forall k \in \mathbb{N}$. Equally, it can be written that $\aleph_1 \cdot (1/k) = \aleph_1$, $\forall k \in \mathbb{N}$. And this last statement can also be extended to the product of $\aleph_1$ by a rational number $q$, i.e. $\aleph_1 \cdot q = \aleph_1$, $\forall q \in \mathbb{Q}: 0 < q < 1$. By considering (4.7) together with Table 4.1, it will then be possible to evaluate the cardinalities of $\mathcal{N}(N_{\mathbb{E},\mathbb{O}+1})$ and subsequent subsets of $\mathcal{P}(\mathbb{N})$:

$$\left|\mathcal{N}(N_{\mathbb{E},\mathbb{O}+1})\right| = \left|\mathcal{N}(N_{\mathbb{E},\mathbb{O}})\right| \cdot q = \lim_{n \to \aleph_0} \left[\binom{n}{n/2} \cdot \frac{1}{1 + (2/n)}\right] = \aleph_1 \tag{4.21}$$

$$\left|\mathcal{N}(N_{\mathbb{E},\mathbb{O}+2})\right| = \left|\mathcal{N}(N_{\mathbb{E},\mathbb{O}+1})\right| \cdot q = \lim_{n \to \aleph_0} \left[\binom{n}{(n/2)+1} \cdot \frac{1 - (2/n)}{1 + (4/n)}\right] = \aleph_1 \tag{4.22}$$

and so on. It is clear that the product of $\aleph_1$ by the corresponding values of the ratio $q$ would render $\aleph_1$ in all cases, since $0 < q < 1$. This is, all the terms in the expansion (4.16) between $\mathcal{N}(N_{\mathbb{E},\mathbb{O}})$ and $\mathcal{N}(N_{\mathbb{N} \setminus p})$ would have cardinalities that equal $\aleph_1$. And, more significantly, the sets $\mathcal{N}(N_{\mathbb{N} \setminus p})$, $\forall p \in \mathbb{N}$, would also have cardinalities equal to $\aleph_1$, therefore contradicting (4.15). Conclusively, such a contradiction confirms that $\mathcal{N}(N_{\mathbb{E},\mathbb{O}})$, and all the subsequents subsets of $\mathcal{P}(\mathbb{N})$, cannot be uncountable.

It is of interest to consider the cardinality of the last term of the expansion (4.16), i.e. $\mathcal{N}(N_{\mathbb{N} \setminus 0})$: If the term $\mathcal{N}(N_{\mathbb{N} \setminus 1})$ had cardinality equal to $\aleph_1$, then the cardinality of $\mathcal{N}(N_{\mathbb{N} \setminus 0})$ would be given by [4]

$$\left|\mathcal{N}(N_{\mathbb{N} \setminus 0})\right| = \left|\mathcal{N}(N_{\mathbb{N} \setminus 1})\right| \cdot q = \lim_{n \to \aleph_0} \left[\aleph_1 \cdot (1/n)\right] = \aleph_1 \cdot (1/\aleph_0) = \aleph_1 \tag{4.23}$$

contradicting the fact that $|\mathcal{N}(N_{\mathbb{N}\setminus 0})| = 1$.

Alternatively, it is easy to see that, if $\mathcal{N}(N_{\mathbb{E},\mathbb{O}})$ had been a denumerable set, i.e. $|\mathcal{N}(N_{\mathbb{E},\mathbb{O}})| = \aleph_0$, then it would had also followed that all the subsequent subsets of $\mathcal{P}(\mathbb{N})$ were denumerable, and that $|\mathcal{N}(N_{\mathbb{N}\setminus 0})| = 1$.

Once all subsets of $\mathcal{P}(\mathbb{N})$, whose members are infinite subsets of $\mathbb{N}$, have been shown to be denumerable, the denumerability of $\mathcal{P}(\mathbb{N})$ follows. Since it can now be written that

$$(4.24) \quad |\mathcal{P}(\mathbb{N})| = 1 + \aleph_0 + \cdots + \aleph_0 + \cdots + \aleph_0 + \cdots + \aleph_0 + \cdots + \aleph_0 + 1$$

where, according to (4.13), the total of terms of the summation is denumerable, a final implementation of Theorem 4.1 concludes that $|\mathcal{P}(\mathbb{N})| = \aleph_0$. □

It becomes a corollary of Theorem 4.2 to state that the set of reals $\mathbb{R}$ is also denumerable [15]:

$$(4.25) \quad |\mathbb{R}| = |\mathcal{P}(\mathbb{N})| = 2^{\aleph_0} = \aleph_0$$

And the same conclusion applies to the power set of $\mathcal{P}(\mathbb{N})$,

$$(4.26) \quad |\mathcal{P}(\mathcal{P}(\mathbb{N}))| = \aleph_0$$

as well as to subsequent power sets, thus questioning the viability of transfinite cardinals beyond $\aleph_0$.

## 5. COMPARING INFINITIES: RELATIVE CARDINALITIES

Since, as a consequence of Theorem 4.2, all common infinite sets appear to have the same cardinality ($\aleph_0$), a different way of comparing them seems necessary. For this purpose, the concept of relative cardinality was introduced in [15], modified here as follows:

**Definition 5.1 (Relative Cardinality of Finite Sets).** *Consider two finite sets, A and B, with cardinalities $|A| = a$ and $|B| = b$, such that $a < b$. Their relative cardinality is defined as the ratio $\rho_{A,B} = a/b$.*

**Definition 5.2 (Relative Cardinality of Infinite Sets).** *Consider two sets, A and B, both denumerable, such that A is a subset of B, $A \subseteq B$. Assume their constructions generate formulae, $\Phi_A(n)$ and $\Phi_B(n)$, $\forall n \in \mathbb{N}$, which render the cardinalities of the respective interim finite sets, in relation to each other.*

---

[4] According to the arithmetics of transfinite cardinals [7,12,16], $\aleph_1 \cdot \aleph_0 = \aleph_1$. This will also allow us to write $\aleph_1 \cdot (1/\aleph_0) = \aleph_1$.

*The relative cardinality of A and B is defined as the limiting ratio*

$$\rho_{A,B} = \lim_{n \to \aleph_0} \frac{\Phi_A(n)}{\Phi_B(n)} \tag{6.1}$$

From both definitions, it will always be the case that $0 \leq \rho_{A,B} \leq 1$. While the definition of relative cardinality for finite sets is elementary, the implementation of (6.1) for infinite sets needs some kind of baseline for $\Phi_A(n)$ and $\Phi_B(n)$ to be truly comparable, hence the requisite of $A$ being a subset of $B$. It can be envisaged that the quality of being a finite or an infinite set is an "absolute" property of the set, based on their definition [7]:

> "A set $X$ is *finite* if there is a bijection $f: n \to X$ for some $n \in \mathbb{N}$. If there is no such bijection for any $n \in \mathbb{N}$, $X$ is *infinite*."

The introduction of relative cardinalities, as defined here, brings the possibility of comparing finite sets according to their size. And the same applies to infinite sets, albeit in relative terms. In this sense, the property of being a denumerable set, i.e. countably infinite, is treated as a basic ("absolute") property of the set which differentiates it from any finite set, but does not from other denumerable sets. However, the relative cardinality $\rho_{A,B}$ of two denumerable sets provides a comparison of their "relative" sizes. This point can better be illustrated by relevant examples:

*i)* Consider the set of natural numbers, $\mathbb{N} = \{0, 1, 2, \cdots, n, \cdots\}$, and the set of even numbers, $\mathbb{E} = \{0, 2, 4, \cdots, 2k, \cdots\} \, \forall k \in \mathbb{N}$. $\mathbb{E}$ is a infinite subset of $\mathbb{N}$, such that $|\mathbb{E}| = |\mathbb{N}| = \aleph_0$, so they have the same "absolute" cardinality. Nevertheless, their relative cardinality is

$$\rho_{\mathbb{E},\mathbb{N}} = \lim_{n \to \aleph_0} \frac{n/2}{n} = 0.5 \tag{6.2}$$

therefore, in relative terms, there are half as many even numbers as natural numbers.

*ii)* Consider the set of natural numbers, $\mathbb{N}$, and the set of all integers, positive and negative, $\mathbb{Z} = \{\cdots, -2, -1, 0, 1, 2, \cdots\}$. $\mathbb{N}$ is a infinite subset of $\mathbb{Z}$, such that $|\mathbb{N}| = |\mathbb{Z}| = \aleph_0$, so once more they have the same "absolute" cardinality. To determine their relative cardinality, $\rho_{\mathbb{N},\mathbb{Z}}$, it is necessary to have a suitable formula $\Phi_{\mathbb{Z}}(n)$. This will be: $\Phi_{\mathbb{Z}}(n) = 2n+1$. Accordingly,

$$\rho_{\mathbb{N},\mathbb{Z}} = \lim_{n \to \aleph_0} \frac{n}{2n+1} = 0.5 \tag{6.3}$$

so, in relative terms, there are half as many natural numbers as integers.

*iii)* Consider the set of natural numbers, $\mathbb{N}$, and the set of rational numbers, $\mathbb{Q} = \{a/b, \forall a, b \in \mathbb{Z}, b \neq 0\}$. $\mathbb{N}$ is a infinite subset of $\mathbb{Q}$, such that $|\mathbb{N}| = |\mathbb{Q}| = \aleph_0$, this is, they have the same "absolute" cardinality. In order to determine their relative cardinality, $\rho_{\mathbb{N},\mathbb{Q}}$, it is necessary to have a suitable formula $\Phi_{\mathbb{Q}}(n)$. Consider first the fractions $q$ in the interval (0,1], such that $0 < q \leq 1, \forall q \in (0,1]$:

(6.4)

| a\b | 1 | 2 | 3 | 4 | 5 | 6 | 7 | 8 | 9 | ... |
|---|---|---|---|---|---|---|---|---|---|---|
| 1 | **1/1** | 2/1 | 3/1 | 4/1 | 5/1 | 6/1 | 7/1 | 8/1 | 9/1 | ... |
| 2 | **1/2** | 2/2 | 3/2 | 4/2 | 5/2 | 6/2 | 7/2 | 8/2 | 9/2 | ... |
| 3 | **1/3** | **2/3** | 3/3 | 4/3 | 5/3 | 6/3 | 7/3 | 8/3 | 9/3 | ... |
| 4 | **1/4** | 2/4 | **3/4** | 4/4 | 5/4 | 6/4 | 7/4 | 8/4 | 9/4 | ... |
| 5 | **1/5** | **2/5** | **3/5** | **4/5** | 5/5 | 6/5 | 7/5 | 8/5 | 9/5 | ... |
| 6 | **1/6** | 2/6 | 3/6 | 4/6 | **5/6** | 6/6 | 7/6 | 8/6 | 9/6 | ... |
| 7 | **1/7** | **2/7** | **3/7** | **4/7** | **5/7** | **6/7** | 7/7 | 8/7 | 9/7 | ... |
| 8 | **1/8** | 2/8 | **3/8** | 4/8 | **5/8** | 6/8 | **7/8** | 8/8 | 9/8 | ... |
| 9 | **1/9** | **2/9** | 3/9 | **4/9** | **5/9** | 6/9 | **7/9** | **8/9** | 9/9 | ... |
| . | . | . | . | . | . | . | . | . | . | ... |
| . | . | . | . | . | . | . | . | . | . | ... |
| . | . | . | . | . | . | . | . | . | . | ... |

(6.4) shows the initial set of fractions $a/b$ for the ranges $1 \leq a,b \leq 9$. Since the elimination of repetitions is required, highlighted in bold are the fractions $a/b$ that have to be accounted for. For the range 1 to $n$, the formula $\Phi_{\mathbb{Q}}(n)$ covering the interval (0,1] will be

(6.5) $$\left(\Phi_{\mathbb{Q}}(n)\right)\Big|_{(0,1]} = [(n^2-n)/2] \cdot f$$

where $f$ represents the average correction factor that discounts repetitions such as 2/4 or 3/6. Figure 2 shows that $f$ takes approximately a value $f = 0.63$. In order to account for the whole of the number line, it is necessary to multiply (6.5) by a factor of $2n$, and add 1 (to account for 0); this is

(6.6) $$\Phi_{\mathbb{Q}}(n) = 2n \cdot \{[(n^2-n)/2] \cdot 0.63\} + 1$$

We are now in a position to evaluate the relativa cardinality $\rho_{\mathbb{N},\mathbb{Q}}$:

(6.7) $$\rho_{\mathbb{N},\mathbb{Q}} = \lim_{n \to \aleph_0} \frac{n}{1.26\,n \cdot [(n^2-n)/2] + 1} = 0$$

This implies that, in relative terms, there are infinitely more rational numbers than natural numbers.

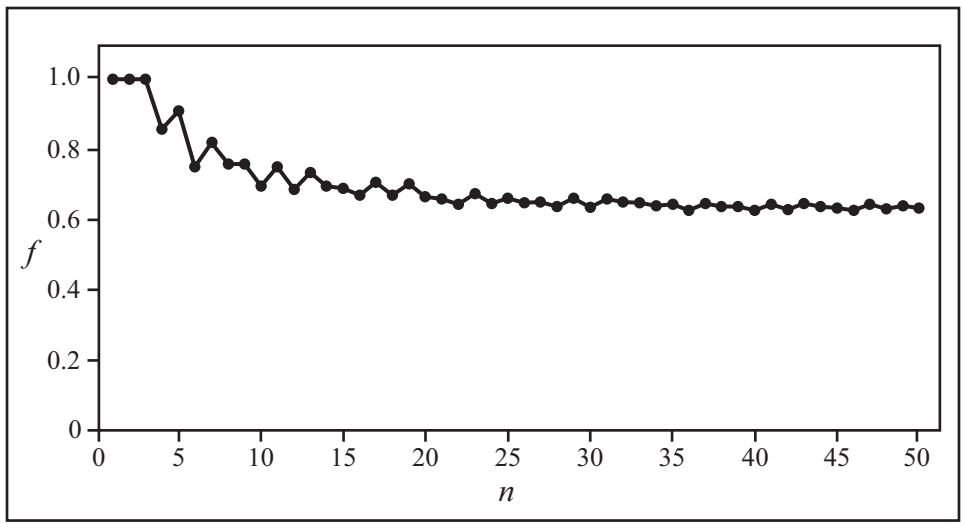

FIGURE 2. Correction factor $f$ as a function of $n$.

*iv)* Consider the set of natural numbers, $\mathbb{N}$, and the set of real numbers, $\mathbb{R}$. $\mathbb{N}$ is a infinite subset of $\mathbb{R}$, such that $|\mathbb{N}|=|\mathbb{R}|=\aleph_0$, this is, they have the same "absolute" cardinality. In order to determine their relative cardinality, $\rho_{\mathbb{N},\mathbb{R}}$, it is first necessary to have a suitable formula $\Phi_{\mathbb{R}}(n)$. Consider the real numbers $r$ in the interval $[0,1)$, such that $0 \leq r < 1, \forall r \in [0,1)$; it will be the case that

$$(6.8) \qquad \left(\Phi_{\mathbb{R}}(n)\right)\Big|_{[0,1)} = 2^n$$

since it is possible to establish a one-to-one correspondence between the set of infinite binary strings and the real numbers in the interval $[0,1)$. To cover the whole of the number line, it will be necessary to multiply (6.8) by $2n$, resulting in

$$(6.9) \qquad \Phi_{\mathbb{R}}(n) = 2n \cdot 2^n = n \cdot 2^{n+1}$$

Finally, the relative cardinality $\rho_{\mathbb{N},\mathbb{R}}$ will be given by

$$(6.10) \qquad \rho_{\mathbb{N},\mathbb{R}} = \lim_{n \to \aleph_0} \frac{n}{n \cdot 2^{n+1}} = 0$$

This implies that, in relative terms, there are infinitely more real numbers than natural numbers.

*iv)* Consider the set of natural numbers, $\mathbb{Q}$, and the set of real numbers, $\mathbb{R}$. $\mathbb{Q}$ is a infinite subset of $\mathbb{R}$, such that $|\mathbb{Q}|=|\mathbb{R}|=\aleph_0$, this is, they have the same "absolute" cardinality. To determine their relative cardinality $\rho_{\mathbb{Q},\mathbb{R}}$, it is sufficient to apply (6.1) using the formulae $\Phi_{\mathbb{Q}}(n)$ and $\Phi_{\mathbb{R}}(n)$ already obtained. The result will be

(6.11) $$\rho_{\mathbb{N},\mathbb{Q}} = \lim_{n \to \aleph_0} \frac{1.26\, n \cdot [(n^2-n)/2] + 1}{n \cdot 2^{n+1}} = 0$$

(6.11) implies the existence, in relative terms, of infinitely more real numbers than rational numbers.

*v)* Consider the set of real numbers, $\mathbb{R}$, and the set of complex numbers, $\mathbb{C}$. $\mathbb{R}$ is a infinite subset of $\mathbb{C}$, such that $|\mathbb{R}| = |\mathbb{C}| = \aleph_0$, this is, they have the same "absolute" cardinality. To determine their relative cardinality $\rho_{\mathbb{R},\mathbb{C}}$, it is first necessary to evaluate the formula $\Phi_{\mathbb{C}}(n)$. Since the set of complex numbers can be constructed by the cartesian product of $\mathbb{R}$ with the set of imaginary numbers $\mathbb{I} = \{r \cdot i, \forall r \in \mathbb{R} \land i = \sqrt{-1}\}$, it can be deduced that

(6.12) $$\Phi_{\mathbb{C}}(n) = (n \cdot 2^{n+1}) \cdot (n \cdot 2^{n+1}) = n^2 \cdot 2^{2n+2}$$

Accordingly, the relative cardinality $\rho_{\mathbb{R},\mathbb{C}}$ will be given by

(6.13) $$\rho_{\mathbb{R},\mathbb{C}} = \lim_{n \to \aleph_0} \frac{n \cdot 2^{n+1}}{n^2 \cdot 2^{2n+2}} = \lim_{n \to \aleph_0} \frac{1}{n \cdot 2^{n+1}} = 0$$

which implies the existence, again in relative terms, of infinitely more complex numbers than real numbers.

Examples *i)* to *v)* illustrate how a concept as simple as the relative cardinality of two given sets $A$ and $B$, $\rho_{A,B}$ (Definitions 6.1 and 6.2) is nevertheless capable of providing a powerful quantitative comparison between the relative sizes of sets, of particular significance when dealing with infinite sets. As Hilbert's well-know metaphor of the "Infinity Hotel" [6,9] indicates, infinity is treated mathematically as an "elastic" entity that can be expanded indefinitely to accommodate more and more members (a property that is fully encapsulated by Theorem 4.1). Since the main claim reported here is that all infinite sets are denumerable, i.e. they all have the same "absolute" cardinality $\aleph_0$, it becomes clear that their relative cardinalities offer effective and quantitative means with which to compare them.

## 6. CONCLUDING REMARKS

The results reported here offer an additional confirmation of the conclusions and implications already reported in [15]. The purge of Cantor's transfinite theory from the fabric of mathematics, although it will undoubtedly be a traumatic and arduous process, will nevertheless bring a great deal of benefits in terms of the consequential simplification of the axiomatic principles that underpin set theory. Such benefits might propagate into all areas of pure mathematics, so only time will tell what new and exciting findings will be uncovered as a result.

## REFERENCES

[1] Andreescu T. and Feng Z. *A Path to Combinatorics for Undergraduates: Counting Strategies*. Boston: Birkhäuser, 2004.
[2] Bartle R.G and Sherbert D.R. *Introduction to Real Analysis*. Third edition. New York: John Wiley & Sons, 2000.
[3] Cantor, G. *Über eine elementare Frage der Mannigfaltigkeitslehre*. Jahresberg. Deutsch. Math.-Ver. 1891; **1**: 75-78.
[4] Coiras E. *Counterexamples to Cantorian Set Theory*. arXiv:1404.6447v1 [maths.GM] 24 April 2014.
[5] Ewald, W. B., ed. *From Kant to Hilbert: A Source Book in the Foundations of Mathematics*. Volume II. Oxford: Clarendon Press, 1996.
[6] Gamov G. *One Two Three... Infinity: Facts and Speculations of Science*. New York: Viking Press, 1947.
[7] Goldrei D. *Classic Set Theory for Guided Independent Study*. Boca Raton: Chapman & Hall/CRC Press, 1998.
[8] Hamilton A.G. *Logic for Mathematicians*. Revised Edition. Cambridge: Cambridge University Press, 1988.
[9] Hilbert, D. *David Hilbert's Lectures on the Foundations of Arithmetics and Logic*. Ewald, W. & Sieg, W. (eds.). Heidelberg: Springer-Verlag, 2013.
[10] Jones, P.P. *The Case Against Cantor's Diagonal Argument*. Research Square: https://doi.org/10.21203/rs.3.rs-40722/v1. 2020.
[11] Knapp C. & Silva C.E. *The Uncountability of the Unit Interval.* arXiv:1209.5199v2 [maths.HO] 22 January 2014.
[12] Lipschutz S. *Set Theory and Related Topics*. Second edition. Schaum's Outline Series. New York: McGraw-Hill, 1998.
[13] Mueckenheim W. *On Cantor's Important Proofs*. arXiv: math/0306200v3 [maths.GM] 17 December 2006.
[14] Nolt J., Rohatyn D. and Varzi A. *Theory and Problems of Logic*. Second edition. Schaum's Outline Series. New York: McGraw-Hill, 1998.
[15] Perez, J.A. *Addressing Mathematical Inconsistency: Cantor and Gödel refuted*. arXiv:1002.433v1 [maths.GM] 25 February 2010.
[16] Potter, M. *Set Theory and its Philosophy*. Oxford: Oxford University Press, 2004.
[17] Smith P. *An Introduction to Gödel's Theorems*. Cambridge: Cambridge University Press, 2007.

March 2023
JUAN A. PEREZ. BERKSHIRE, UK.
*Email address*: jap717@juanperezmaths.com